\documentclass[a4paper,12pt, reqno]{amsart}
\usepackage{amsthm}
\usepackage{amssymb}

\newtheorem{thm}{Theorem}[section]

\newtheorem{lem}[thm]{Lemma}

\theoremstyle{definition}

\theoremstyle{remark}

\begin{document}
%%% \topmatter
\title[Non-Wieferich Prime]{Non-Wieferich primes in number fields and \emph{ABC} conjecture}
\author{Srinivas Kotyada and Subramani Muthukrishnan}
\address{Institute of Mathematical Sciences, HBNI, 
CIT Campus, Taramani, Chennai 600 113, India}
\address{Chennai Mathematical Institute,
SIPCOT IT Park, Siruseri, Chennai 603 103, India}
\email[Kotyada Srinivas]{srini@imsc.res.in}
\email[Subramani Muthukrishnan]{subramani@cmi.ac.in}
\begin{abstract}
Let $K/\mathbb{Q}$ be an algebraic number field of  class number one and  $\mathcal{O}_K$ be its ring of integers.  We show that there are infinitely many non-Wieferich primes  with respect to certain units in $\mathcal{O}_K$ under the assumption of the $\emph{abc}$ conjecture for number fields.
\end{abstract}

\subjclass[2010]{11A41 (primary); 11R04
(secondary)}
\maketitle

\section{Introduction}

An odd rational prime $p$ is called Wieferich prime if 
\begin{equation}
2^{p-1} \equiv 1 \pmod {p^2}.
\end{equation}
A. Wieferich \cite{Wie} proved that if an odd prime $p$ is non-Wieferich prime, i.e., $p$  satisfies 
$$ 
2^{p-1} \not \equiv 1 \pmod {p^2},
$$ 
then there are no integer solutions to the Fermat equation $x^p+y^p=z^p,$ with $p \nmid xyz$. The known Wieferich primes are $1093$ and $3511$ and  according to the PrimeGrid project \cite{PrimeGrid}, these are the only Wieferich primes less than $17 \times 10^{15}$. One of the unsolved problems in this area of research is to determine whether the number of Wieferich or non-Wieferich primes is finite or infinite. Instead of the base $2$ if we take any base $a$, then $p$ is said to be a Wieferich prime with respect to the base $a$ if
\begin{equation}\label{non-w}
a^{p-1} \equiv 1 \pmod {p^2},
\end{equation}
and if the congruence \eqref{non-w} does not hold then we shall say that $p$ is non-Wieferich prime to the base $a$. Under the famous \emph{abc} conjecture (defined below), J. H. Silverman \cite{Sil} proved that given any integer $a$, there are infinitely many non-Wieferich primes to the base $a$.  He established this result by showing that for any fixed $\alpha \in \mathbb{Q}^{\times},  \alpha \not = \pm 1,$ and assuming the truth of $\emph{abc}$ conjecture, 
$$ \mbox{card} \left\{p \leq x : \alpha^{p-1} \not \equiv 1 \pmod{p^2} \right\} \gg_{\alpha} \log x \quad \textrm{as} \quad x \to \infty.$$ 
In \cite{Graves} Hester Graves and M. Ram Murty extended this result to primes in arithmetical progression by showing that  for any $a \geq 2$ and any fixed $k \geq 2$, there are $\gg \log x/{\log\log x}$ primes $p\leq x$ such that $a^{p-1} \not \equiv 1 \pmod{p^2}$ and $p\equiv 1 \pmod{k}$, under the assumption of $\emph{abc}$ conjecture.

\medskip

\noindent
In this paper, we study non-Wieferich primes in algebraic number fields of class number one.
More precisely, we prove the following 
\begin{thm}\label{thm-1}
Let $K=\mathbb{Q}(\sqrt{m})$ be a real quadratic field of class number one and assume that the $\emph{abc}$ conjecture holds true in $K$. Then there are infinitely many non-Wieferich primes in $\mathcal{O}_K$ with respect to the unit $\varepsilon$ satisfying $|\varepsilon|>1$.
\end{thm}
\noindent

\begin{thm}\label{thm-2}
Let $K$ be any algebraic number field of class number one and assume that the $\emph{abc}$ conjecture holds true in $K$. Let $\eta$ be a unit in $\mathcal{O}_K$ satisfying $|\eta|>1$ and $|\eta^{(j)}|<1$ for all $j \not = 1$, where $\eta^{(j)}$ is the $j$th conjugate of $\eta$. Then there exist infinitely many non-Wieferich primes in $K$ with respect to the base $\eta$.
\end{thm}

\medskip

\noindent
The plan of this article is as follows. In section 2, we shall define the $abc$ conjecture for number fields. In section 3, a brief introduction to Wieferich/non-Wieferich primes over number fields will be given and in section 4 and 5, we shall prove Theorem 1.1 and Theorem 1.2, respectively.

\medskip

\section{The $abc$-conjecture}

\noindent
The ${abc}$-conjecture propounded by Oester\'e and Masser (1985) states that  given any $\delta >0$ and positive integers $a,b,c$ such that $a+b=c$ with $(a,b)=1$, we have
$$ c\ll_{\delta} (\mbox{rad}(abc))^{1+\delta},$$
where $\mbox{rad} (abc) := \prod_{p|abc}p$.

\medskip

\noindent
The $abc$ conjecture has several applications,  the reader may refer to \cite{Vojta}, \cite{Gyory}, \cite{ABC}, \cite{ABC-Ram}  for details. 

\medskip

\noindent
To state the analogue of $abc$-conjecture for number fields, we need some preparations, which we do below. The interested reader may refer to \cite{Vojta}, \cite{Gyory} for more details.
\medskip

\noindent
Let $K$ be an algebraic number field and let $V_{K}$ denote the set of primes on $K$, that is, any $v$ in $V_K$ is an equivalence class of norm on $K$ (finite or infinite). Let $||x||_v := N_{K/\mathbb{Q}}(\mathfrak{p})^{-v_\mathfrak{p}(x)},$ if $v$ is a prime defined by the prime ideal $\mathfrak{p}$ of the ring of integers $\mathcal{O}_K$ in $K$ and $v_\mathfrak{p}$ is the corresponding valuation, where $N_{K/\mathbb{Q}} $ is the absolute value norm. Let $||x||_v := |g(x)|^e$ for all non-conjugate embeddings $ g: K \to \mathbb{C}$ with $e=1$ if g is real and $e=2$ if g is complex. Define the height of any triple $a,b,c \in K^{\times}$ as 
\begin{equation*}
H_K(a,b,c) := \prod _{v \in V_K} max (||a||_v, ||b||_v, ||c||_v),
 \end{equation*}

\noindent
and the radical of $(a,b,c)$ by
$$\mbox{rad}_K(a,b,c) := \prod_{\mathfrak{p} \; \in \; I_K(a,b,c)}	N_{K/\mathbb{Q}}(\mathfrak{p})^{v_\mathfrak{p}{(p)}},
$$
where $p$ is a rational prime with $p\mathbb{Z}= \mathfrak{p} \cap \mathbb{Z}$ and   $I_K(a,b,c)$ is the set of all primes $\mathfrak{p}$ of $\mathcal{O}_K$ for which $||a||_v, ||b||_v, ||c||_v$ are not equal. 

\medskip

\noindent
The abc conjecture for algebraic number fields is stated as follows: For any $\delta >0$, we have
\begin{equation}\label{abc-for-nf}
H_K(a,b,c) \ll_{\delta,K} (\mbox{rad}_K (a,b,c))^{1+\delta},
\end{equation}
for all $a,b,c, \in K^{\times}$ satisfying $a+b+c=0$, the implied constant depends on $K$ and $\delta$.

\section{Wieferich/non-Wieferich primes in number fields}

\noindent
Let $K$ be an algebraic number field and $\mathcal{O}_K$ be its ring of integers. A prime $\pi \in \mathcal{O}_K$ is called Wieferich prime with respect to the base $\varepsilon \in \mathcal{O}_K^{*}$ if
\begin{equation}\label{def-w}
\varepsilon^{N(\pi)-1} \equiv 1 \pmod{\pi^2},
\end{equation}
where $N(.)$ is the absolute value norm.  If the congruence \eqref{def-w} does not hold for a prime $\pi \in \mathcal{O}_K$, then it is  called non-Wieferich prime to the base $\varepsilon.$

\medskip

\noindent
\textbf{Notation:} In what follows,  $\varepsilon$ will denote a unit in $\mathcal{O}_K$ and we shall write $\varepsilon^n - 1 = u_n v_n$, where $u_n$ is the square free part and $v_n$ is the squarefull part, i.e.,  if $\pi | v_n$ then $\pi^2 | v_n$. We shall denote absolute value norm on $K$ by $N$.

\section{Proof of theorem \eqref{thm-1}}

\noindent
Let $K=\mathbf{Q}(\sqrt{m}), m>0$ be a real quadratic field and $\mathcal{O}_{K}$ be its ring of integers. Let $\varepsilon \in \mathcal{O}_K^*$ be a unit with $|\varepsilon| >1$. The results of Silverman \cite{Sil}, Ram Murty and Hester \cite{Graves} elucidated in the introduction use a key lemma of Silverman (Lemma 3, \cite{Sil}). We derive an analogue of Silverman's lemma for number fields which will play a fundamental role in the proof of the main theorems.

\begin{lem}\label{lemma-1}
Let $K = \mathbf{Q}(\sqrt{m})$ be a real quadratic field of class number one. Let $\varepsilon \in \mathcal{O}_K^*$ be a unit. If $\varepsilon^{n}-1 = u_n v_n,$ then every prime divisor $\pi$ of $u_n$ is a non-Wieferich prime with respect to the base $\varepsilon$.
\end{lem}

\noindent
\textbf{Proof.}  The assumption that $K$ has class number one allows us to write the element $\varepsilon^{n}-1 \in \mathcal{O}_K$ as a product of primes uniquely. Accordingly, we shall write
 $$\varepsilon^{n}-1 = u_n v_n$$ for $n \in \mathbb{N}$. Then  
\begin{equation}\label{eq-1}
\varepsilon^n = 1 + \pi w,
\end{equation}
with  $\pi | u_n$ and $\pi$ and $w$ are coprime.
As $\pi$ is a prime, we have $N(\pi) = p$ or $p^2$, $p$ is a rational prime.

\medskip

\noindent
Case (1): Suppose $N(\pi) = p$.

\noindent
From equation \eqref{eq-1},  we get 
$$
\varepsilon^{n(p-1)} \equiv 1+(p-1) \pi w \not \equiv 1 \pmod{\pi^2}.
$$
\noindent
Case (2): Suppose $N(\pi) = p^2.$

\noindent
Again from  equation \eqref{eq-1}, we obtain
$$ 
{\varepsilon^n}^{(p^2-1)} = {\varepsilon^n}^{(N(\pi)-1)} = (1+\pi w)^{(p^2-1)} \equiv 1 +\pi w (p^2-1) \not \equiv 1 \pmod{\pi^2}.
$$
\noindent
Thus in either case, 
$$
\varepsilon^{(N(\pi)-1)} \not \equiv 1 \pmod{\pi^2},
$$ 
and hence $\pi$ is a non-Wieferich prime to the base $\varepsilon.$ \qed

\medskip

\noindent
The above lemma shows that whenever a prime $\pi$ divides $u_n$ for some positive integer $n$, then $\pi$ is a non-Wieferich prime with respect to the base $\varepsilon$. Thus, if we can show that the set $\{N(u_n) : n\in \mathbb{N}\}$ is unbounded, then this will imply that the set $\{\pi : \pi | u_n, n\in \mathbb{N} \}$ is an infinite set. Consequently, this establishes the fact that there are infinitely many non-Wieferich primes in every real quadratic field of class number one with respect to the unit $\varepsilon$, with $| \varepsilon | > 1$. Therefore, we need only to show the following

\begin{lem}
Let $\mathbb{Q}(\sqrt{m})$ be a real quadratic field of class number one. Let $\varepsilon \in \mathcal{O}_K^*$ be a unit with $|\varepsilon|>1$. Then under $abc$-conjecture for number fields, the set $\{N(u_n) : n\in \mathbb{N}\}$ is unbounded.
\end{lem}
\noindent
\textbf{Proof.} Invoking  the $abc$-conjecture \eqref{abc-for-nf}  to the equation
\begin{equation}\label{basic-eq}
\varepsilon^n = 1 + u_n v_n
\end{equation}
yields
\begin{equation}\label{basic-2}
|\varepsilon^n| \ll \big(\prod_{\mathfrak{p}|u_n v_n} N(\mathfrak{p})^{v_\mathfrak{p}(p)}\big)^{1+\delta} = \bigg(\prod_{\mathfrak{p}|u_n} N(\mathfrak{p})^{v_\mathfrak{p}(p)} \prod_{\mathfrak{p}|v_n} N(\mathfrak{p})^{v_\mathfrak{p}(p)} \bigg)^{1+\delta}
\end{equation}
for some $\delta > 0$. Here the implied constant depends on $K$ and $\delta$. 

\medskip
\noindent
As $v_\mathfrak{p}(p)\leq 2$ for any prime ideal $\mathfrak{p}$ lying above the rational prime $p,$ we have

\begin{equation}\label{neweqn-1}
\prod_{\mathfrak{p}|u_n} N(\mathfrak{p})^{v_\mathfrak{p}(p)} \leq N(u_n)^2.
\end{equation}

\medskip
\noindent
For a prime ideal $\mathfrak{p}|v_n,$ let $e_\mathfrak{p}$ be the largest exponent of $\mathfrak{p}$ dividing $v_n,$ i.e., $\mathfrak{p}^{e_\mathfrak{p}}||v_n.$ As $v_n$ is the square-full part of $\varepsilon^n-1,$ we have $e_\mathfrak{p} \geq 2$. Hence,

\medskip

\begin{enumerate}
\item $N(\mathfrak{p})^{2v_\mathfrak{p}(p)} \leq N(\mathfrak{p})^{2+e_\mathfrak{p}}$ for all  prime ideals $\mathfrak{p}$ with $v_\mathfrak{p}(p)=2.$ 
\item  $N(\mathfrak{p})^{2v_\mathfrak{p}(p)} \leq N(\mathfrak{p})^{e_\mathfrak{p}}$ for all prime ideals $\mathfrak{p}$ with $v_\mathfrak{p}(p)=1.$ 
\end{enumerate}
Thus

\begin{align*}
\prod_{\mathfrak{p}|v_n} N(\mathfrak{p})^{2v_\mathfrak{p}(p)} & \leq  \prod_{\mathfrak{p}|v_n \atop v_\mathfrak{p}(p)=2} N(\mathfrak{p})^{2+e_\mathfrak{p}(p)} \prod_{\mathfrak{p}|v_n \atop v_\mathfrak{p}(p)=1} N(\mathfrak{p})^{e_\mathfrak{p}(p)}\\
& \leq  \prod_{\mathfrak{p}|v_n \atop v_\mathfrak{p}(p)=2} N(\mathfrak{p})^2 \prod_{\mathfrak{p}|v_n \atop v_\mathfrak{p}(p)=2} N(\mathfrak{p})^{e_\mathfrak{p}(p)}\prod_{\mathfrak{p}|v_n \atop v_\mathfrak{p}(p)=1} N(\mathfrak{p})^{e_\mathfrak{p}(p)}\\
&\leq  {\prod_{\mathfrak{p}}}' N(\mathfrak{p})^2 \prod_{\mathfrak{p}|v_n \atop v_\mathfrak{p}(p)=2} N(\mathfrak{p})^{e_\mathfrak{p}(p)}\prod_{\mathfrak{p}|v_n \atop v_\mathfrak{p}(p)=1} N(\mathfrak{p})^{e_\mathfrak{p}(p)},
\end{align*}

\noindent
where $\mathbb{'}$ indicates that the product is over all primes $\mathfrak{p}$ in $\mathcal{O}_K$ such that $v_\mathfrak{p}(p)=2$. As it is well known that there are only finitely many ramified primes in a number field, it follows that the product is bounded by a constant $A$  (say). Thus, we have 

\begin{equation}\label{neweqn-2}
\prod_{\mathfrak{p}|v_n} N(\mathfrak{p})^{v_\mathfrak{p}(p)} \leq \sqrt{AN(v_n)}.
\end{equation}

\medskip
\noindent
Combining equations \eqref{basic-2}, \eqref{neweqn-1} and \eqref{neweqn-2}, we get

\begin{equation}\label{neweqn-3}
|\varepsilon^n| \ll \left(N(u_n)^2 \sqrt{N(v_n)}\right)^{1+\delta}. 
\end{equation}

\medskip

\noindent
Now, as $|\varepsilon| >1$,  

$$N(u_n)N(v_n) = N(\varepsilon^n -1)  \leq 2 |\varepsilon^n-1| < 2 |\varepsilon|^n,$$
i.e.,
$$N(v_n) < 2 |\varepsilon|^n/ N(u_n).$$
Substituting the above expression in \eqref{neweqn-3}, we obtain
$$|\varepsilon^n| \ll \bigg(N(u_n)^2 \frac{|\varepsilon|^{n/2}}{\sqrt{N(u_n)}}\bigg)^{(1+\delta)}. $$
Thus,
$$ (N(u_n))^{\frac{3(1+\delta)}{2}} \gg |\varepsilon|^{\frac{ n(1-\delta)}{2}}.$$
Thus, for a fixed $\delta$, $N(u_n) \to \infty$ as $ n \to \infty$.
This proves the lemma and hence completes the proof of the theorem.\qed

\section{Non-Wieferich primes in algebraic number fields}

In this section, we generalize the arguments of previous section to arbitrary number fields. From now onwards,  $K$ will always denote an algebraic number field of degree $[K:\mathbb{Q}]=l$ over $\mathbb{Q}$ of class number one. Let $r_1$ and $r_2$ be the number of real and non-conjugate complex embeddings of $K$ into $\mathbb{C}$ respectively, so that $ l = r_1 + 2r_2$.  We begin with an analogue of Lemma \eqref{lemma-1}.
\begin{lem}\label{lemma-3}
Let $\varepsilon$ be a unit in $\mathcal{O}_K$. If $\varepsilon^{n}-1 = u_n v_n,$ then every prime divisor $\pi$ of $u_n$ is a non-Wieferich prime with respect to the base $\varepsilon$.
\end{lem}

\noindent
\textbf{Proof.} Let $N(\pi) = p^k,$ where $p$ is a rational prime and $k$ is a positive integer. Then
$$\varepsilon^{n(N(\pi)-1)} = \varepsilon^{n(p^k-1)} = (1+w \pi)^{(p^k-1)} \equiv 1 + (p^k-1) w \pi \not \equiv 1 \pmod{\pi^2}.$$
This implies $\varepsilon^{N(\pi) -1} \not \equiv \: 1 \pmod{\pi^2}.$ 

\medskip

\noindent
Thus, the lemma shows that $\pi$ is a non-Wieferich prime to the base $\varepsilon$ whenever the hypothesis of the lemma is met. Now, under the $abc$ conjecture for number fields, we show below the existence of infinitely many non-Wieferich primes.

\begin{lem}\label{lemma-4}
The set $\{N(u_{n}):n\in \mathbb{N}\}$ is unbounded, where ${u_n}$'s are as defined in Lemma \eqref{lemma-3}.
\end{lem}

\noindent
\textbf{Proof.}  By the hypothesis of the lemma, we have $\varepsilon^n =1 + u_n v_n$, where $\varepsilon^n, 1, u_n v_n \in K^{\times}.$
Applying the $\emph{abc}$ conjecture for number fields to the above equation, we obtain
\begin{equation}\label{eqn-5}
\prod_{v \in V_{K}} \max (|u_n v_n|_v, |1|_v, |\varepsilon^n|_v) \ll (\prod_{\mathfrak{p}|u_n v_n} N(\mathfrak{p})^{v_\mathfrak{p}(p)})^{1+\delta},
\end{equation}
for some $\delta >0.$

\noindent
Note that for the absolute value $|.|$ in $ V_K,$ we have
\begin{equation}\label{eqn-6}
|\varepsilon^n| \leq \prod_{v \in V_{K}} \max (|u_n v_n|_v, |1|_v, |\varepsilon^n|_v).
\end{equation}

\noindent
As $v_\mathfrak{p}(p) \leq l$ for any prime ideal $\mathfrak{p}$ lying above the rational prime $p,$ we have
\begin{equation}
\prod_{\mathfrak{p}|u_n} N(\mathfrak{p})^{v_\mathfrak{p}(p)} \leq N(u_n)^{l}.
\end{equation}

\noindent
As before, we denote by $e_\mathfrak{p}$ the largest exponent of $\mathfrak{p}$ which divides $v_n,$ i.e., $\mathfrak{p}^{e_\mathfrak{p}} || v_n.$ Clearly $e_{\mathfrak{p}} \geq 2.$ Then

\begin{align*}
\prod_{\mathfrak{p}|v_n} N(\mathfrak{p})^{2v_\mathfrak{p}(p)} 
& \leq \prod_{\mathfrak{p}|v_n \atop v_\mathfrak{p}(p) \geq 2} N(\mathfrak{p})^{2l+e_\mathfrak{p}(p)} \prod_{\mathfrak{p}|v_n \atop v_\mathfrak{p}(p)=1} N(\mathfrak{p})^{e_\mathfrak{p}(p)}\\
&\leq \prod_{\mathfrak{p}|v_n \atop v_\mathfrak{p}(p) \geq 2} N(\mathfrak{p})^{2l} \prod_{\mathfrak{p}|v_n \atop v_\mathfrak{p}(p) \geq 2} N(\mathfrak{p})^{e_\mathfrak{p}(p)}\prod_{\mathfrak{p}|v_n \atop v_\mathfrak{p}(p) = 1} N(\mathfrak{p})^{e_\mathfrak{p}(p)}\\
&\leq {\prod_{\mathfrak{p}}}' N(\mathfrak{p})^{2l} \prod_{\mathfrak{p}|v_n \atop v_\mathfrak{p}(p) \geq 2} N(\mathfrak{p})^{e_\mathfrak{p}(p)}\prod_{\mathfrak{p}|v_n \atop v_\mathfrak{p}(p) = 1} N(\mathfrak{p})^{e_\mathfrak{p}(p)},
\end{align*}
where $\mathbb{'}$ indicates that the product is over all primes $\mathfrak{p}$ in $\mathcal{O}_K$ such that $v_\mathfrak{p}(p) \geq 2$. As there are only finitely many ramified primes in a number field, it is bounded by a constant $B$ (say). Thus,  we have

\begin{equation}\label{eqn-7}
\prod_{\mathfrak{p}|v_n} N(\mathfrak{p})^{v_\mathfrak{p}(p)} \leq \sqrt{BN(v_n)}.
\end{equation}

\noindent
Therefore, the equations \eqref{eqn-5} - \eqref{eqn-7} yield
\begin{equation}\label{eqn-8}
|\varepsilon^n| \ll \left(N(u_n)^l\sqrt{N(v_n)}\right)^{1+\delta}.
\end{equation}

\noindent
Note that in the case of real quadratic fields, the unit $\varepsilon$ satisfies $ | \varepsilon | > 1$ and this information was crucial in proving Theorem \ref{thm-1}. However, in the case of general number fields, the following result (see Lemma 8.1.5, \cite{Ram}) comes to our rescue. We state this result as 
\begin{lem} \label{lemma-Ram}
Let $E = \{k\in \mathbb{Z}: 1\leq k \leq r_1+r_2\}.$ Let $E=A \cup B$ be a proper partition of $E.$ There exists a unit $\eta \in \mathcal{O}_K$ with $|\eta^{(k)}| < 1,$ for $k \in A$ and $|\eta^{(k)}|> 1,$ for $k \in B.$
\end{lem}
\noindent
Taking $A = \{k : 1 < k \leq r_1 + r_2 \}$ and $B = \{ 1 \}$,  Lemma \ref{lemma-Ram} produces a unit $\eta \in \mathcal{O}_K^{*}$ such that
$ | \eta | > 1 $ and $ | \eta^{(k)} | <1$, where $\eta^{(k)}$ denotes the $k^{\mathrm{th}}$ conjugate of $\eta, k\neq 1$. 
Since, every unit satisfies \eqref{eqn-8}, replacing $\varepsilon$ with $\eta$ in \eqref{eqn-8}, we obtain
\begin{equation}\label{eqn-9}
|\eta^n| \ll \left(N(u_n)^l\sqrt{N(v_n)}\right)^{1+\delta},
\end{equation}
where, by abuse of notation, we shall denote $\eta^n - 1 = u_n v_n$, with $u_n$ and $v_n$ denoting the same quantities as defined earlier.

\medskip

\noindent
Now,
$$N(u_{n})N(v_{n}) = N(\eta^{n}-1) = (\eta^{n}-1)(\eta^{(2)n}-1)(\eta^{(3)n}-1)\cdots (\eta^{(l)n}-1).
$$
By Lemma \ref{lemma-Ram}, $| \eta^{(j)n}-1| < 2$ for all $j,  2 \leq j \leq l$.

\medskip

\noindent
Thus, 
$$N(u_{n})N(v_{n}) < C |\eta^{n}| \qquad \mathrm{or}  \qquad N(v_{n}) <  C |\eta^{n}|/N(u_{n}).
$$

\noindent
Now, \eqref{eqn-9} can be written as
\begin{equation}\label{last}
\left(N(u_{n})\right)^{\frac{(2l-1)(1+\delta)}{2}} \gg |\eta |^{n\frac{1-\delta}{2}}.
\end{equation}

\noindent
For a fixed $\delta$, the right hand side of \eqref{last} tends to $\infty$ as $ n \to \infty$. 
Therefore the set $\{N(u_{n}):n\in \mathbb{N}\}$ is unbounded.  This shows that there are infinitely many non-Wieferich primes in $K$ with respect to the base $\eta$.  \qed

\medskip

\noindent
\textbf{Acknowledgement}

\noindent
We express our indebtedness to Prof. M. Ram Murty for initiating us into this project and for having many fruitful discussions.  The second author would like to thank Prof. T.R. Ramadas for encouragement and also acknowledges him with thanks for the financial support extended by DST through his J.C. Bose Fellowship. Our sincere thanks to the referee for pointing out some errors and suggesting some changes in an earlier version of this paper. This work is part of the Ph D thesis of the second author.

\end{document}